\documentclass[12pt, a4paper]{amsart}

\usepackage{graphicx}

\usepackage{biograph}      

\usepackage{amssymb}

\usepackage{amsmath}

\usepackage{mathrsfs}      

\usepackage{mathptmx}      

\usepackage[all]{xy}
\CompileMatrices

\newtheorem{theorem}{Theorem}[section]
\newtheorem{lemma}[theorem]{Lemma}
\newtheorem{corollary}[theorem]{Corollary}
\newtheorem{proposition}[theorem]{Proposition}
\newtheorem{definition}[theorem]{Definition}

\newtheorem{example}[theorem]{Example}
\newtheorem{problem}[theorem]{Problem}
\newtheorem{remark}[theorem]{Remark}

\numberwithin{figure}{section}

\DeclareMathOperator{\aut}{Aut}
\DeclareMathOperator{\inn}{Inn}
\DeclareMathOperator{\QC}{QC}
\DeclareMathOperator{\Th}{Th}

\DeclareMathOperator{\autfL}{Autf_L}
\DeclareMathOperator{\galL}{Gal_L}

\DeclareMathOperator{\autfKP}{Autf_{KP}}

\DeclareMathOperator{\autfSh}{Autf_{Sh}}

\DeclareMathOperator{\stab}{Stab}

\DeclareMathOperator{\tp}{tp}
\DeclareMathOperator{\diam}{diam}
\DeclareMathOperator{\lin}{lin}
\DeclareMathOperator{\inter}{int}

\DeclareMathOperator{\acleq}{acl^{eq}}

\DeclareMathOperator{\acl}{acl}

\DeclareMathOperator{\cl}{cl}
\DeclareMathOperator{\id}{id}

\newcommand{\C}{{\mathfrak C}}
\newcommand{\M}{{\mathcal M}}
\newcommand{\F}{{\mathbb F}}
\newcommand{\T}{{\mathcal T}}
\newcommand{\CC}{{\mathcal C}}
\newcommand{\G}{{\mathcal G}}
\newcommand{\B}{{\mathcal B}}

\newcommand{\CCC}{{\mathbb C}}

\newcommand{\V}{{\mathcal V}}

\begin{document}

\title{$G$-compactness and groups}
\author{Jakub Gismatullin \and Ludomir Newelski}
\thanks{The first author is supported by the Polish Goverment grant N N201 384134. The second author is supported by the Polish Goverment grant N201 032 32/2231}
\address{Instytut Matematyczny, Uniwersytetu Wroc\l awskiego,\\ pl. Grunwaldzki 2/4, 50-384 Wroc\l aw, Poland}
\email{gismat@math.uni.wroc.pl \and newelski@math.uni.wroc.pl}


\begin{abstract}
   Lascar described $E_{KP}$ as a composition of $E_{L}$ and the topological
closure of $E_{L}$ (\cite{pillay}). We generalize this result to some other pairs of
equivalence relations.

   Motivated by an attempt to construct a new example of a non-$G$-compact
theory, we consider the following example.
   Assume $G$ is a group definable in a structure $M$. We define a structure
$M'$ consisting of $M$ and $X$ as two sorts, where $X$ is an affine copy of $G$
and in $M'$ we have the structure of $M$ and the action of $G$ on $X$.
   We prove that the Lascar group of $M'$ is a semi-direct product of the
Lascar group of $M$ and $G/G_L$. We discuss the relationship between
$G$-compactness of $M$ and $M'$. This example may yield new examples of
non-$G$-compact theories.
\end{abstract}

\maketitle

\section{Introduction}

Let $T$ be a complete theory in language $L$. We work within a monster model $\C \models T$. A model $M\models T$ is small if $M\prec\C$ and $|M|=|T|$. If $X$ is a subset of a topological space, then by $\inter(X)$ we denote its interior and by $\cl(X)$ its closure. We recall some well known facts about the Lascar Group and Lascar strong types (see \cite{pillay,ziegler}).
The group of Lascar strong automorphisms is defined by: 
\[  \autfL(\C) = \big\langle \aut(\C/M) : M\textrm{ is a model } \big\rangle, \]
and the \emph{Lascar (Galois) group} of $T$ by:
\[ \galL(T) = \aut(\C)/ \autfL(\C).  \]
This definition does not depend on the choice of the monster model $\C$ of $T$ (it is enough that $\C$ is $|T|^+$-saturated and $|T|^+$-strongly homogeneous). We say that $a,b\in \C^k$ ($k<|T|^+$) have the same Lascar strong type, and write $E_L(a,b)$, if there exists $f\in\autfL(\C)$ such that $a=f(b)$. Thus $E_L$ is a $\emptyset$-invariant and bounded equivalence relation on every sort $\C^k$ (because if $a\underset{M}{\equiv}b$ for some small $M\prec \C$, then $E_L(a,b)$, so $|\C^k/E_L| \le |S_k(M)|\le2^{|T|}$).

\begin{definition} \label{def:thick}
A symmetric formula $\varphi(x,y)\in L_{k+k}(\emptyset)$ is \emph{thick} if for some $n<\omega$, for every sequence $(a_i)_{i<n}$ there exist $i<j<n$ such that $\varphi(a_i,a_j)$.  By $\Theta$ we denote the conjunction of all thick formulas:
\[\Theta(x,y) = \bigwedge_{\varphi\textrm{ thick}}\varphi(x,y).\]
\end{definition}

In the above definition we can equivalently take an infinite sequence $(a_i)_{i<\omega}$. If $\varphi$ and $\theta$ are thick, then $\psi(x,y)=\varphi(y,x)$ and $\varphi \land \theta$ are also thick (this follows from the Ramsey theorem). $\Theta$ is a $\emptyset$-invariant relation (not necessarily an equivalence relation) and if $\Theta(a_0,a_1)$, then we can extend $(a_0,a_1)$ to an order indiscernible sequence $(a_i)_{i<\omega}$. On the other hand if $(a_i)_{i<\omega}$ is a 2-indiscernible sequence, then $\varphi(a_0,a_1)$ for every thick $\varphi$.

\begin{lemma}\cite[Lemma 7]{ziegler} \label{lem:zieg}
\begin{enumerate} 
\item[(i)] If $\Theta(a,b)$, then there is a small $M$ such that $a\underset{M}{\equiv}b$.
\item[(ii)] If for some small $M$ we have $a\underset{M}{\equiv}b$, then $\Theta^2(a,b)$, i.e. there is $c$ such that $\Theta(a,c)\land\Theta(c,b)$.
\item[(iii)] $E_L$ is the transitive closure of $\Theta$.
\end{enumerate}
\end{lemma}

If $\pi$ is a type over $\emptyset$, then we can define thick formulas on $\pi(\C)$ and their conjunction $\Theta_{\pi}$ similarly as in the above definition. Moreover, the last remark also holds for $\Theta_{\pi}$, so ${E_L}|_{\pi(\C)}$ is the transitive closure of $\Theta_{\pi}$. One can prove that ${E_L}|_{\pi(\C)}$ is the finest bounded $\emptyset$-invariant equivalence relation on $\pi(\C)$.

There is a compact (not necessarily Hausdorff) topology on the group $\galL(T)$. Let $M$ and $N$ be arbitrary small models and let \[ S_M(N)=\{\tp(M'/N) : \tp(M')=\tp(M)\} \] be a closed subset of $S_{|T|}(N)$. Thus $S_M(N)$ carries a compact subspace topology. The quotient map $j\colon \aut(\C)\rightarrow\galL(T)$ factors as $j=\nu\circ\mu$, where $\mu\colon\aut(\C)\rightarrow S_M(N)$ maps $f$ to $\tp(f(M)/N)$, and $\mu\colon S_M(N)\rightarrow \galL(T)$ maps $\tp(f(M)/N)$ to an appropriate coset of $\autfL(\C)$, so we have the following commutative diagram:
\[\xymatrix{
\aut(\C)  \ar@{>>}[rr]^-{j=\nu\circ\mu} \ar@{>>}[rdd]_-{\mu} &    &  \galL(T) = \aut(\C)/ \autfL(\C) \\
 \\
 & S_M(N) \ar@{>>}[ruu]_-{\nu}
}\]
We can induce topology on $\galL(T)$ from $\nu$, i.e. $X\subseteq\galL(T)$ is closed if and only if its preimage $\nu^{-1}[X]$ is closed in $S_M(N)$. It can be easily seen that this definition of topology does not depend on the choice of small models $M$ and $N$ (\cite[Theorem 4]{ziegler}). With this topology $\galL(T)$ becomes a compact topological group. We say that $T$ is $G$-compact when $\galL(T)$ is Hausdorff. If we consider $\aut(\C)$ with the usual topology of pointwise convergence, then all the maps in the diagram are continuous. However $\nu$ need not be open, instead $\nu$ satisfies some weak kind of openness.

\begin{theorem}\cite[Lemma 12]{ziegler} \label{th:ziegler}
For $p\in S_M(N)$ define its $\Theta$-neighbourhood as:
\[[p]_{\Theta}=\{q\in S_M(N) : p(x)\cup q(y) \cup \Theta(x,y) \textrm{ is consistent }\}. \]
If we take an arbitrary point $p\in S_M(N)$ and subset $U\subseteq S_M(N)$ such that $[p]_{\Theta}\subseteq\inter(U)$, then $\nu(p)\in \inter(\nu[U])$.
\end{theorem}

The relation $E_L$ is $\emptyset$-invariant, so we may consider $E_L$ as a subset of $S_{|T|+|T|}(\emptyset)$. Using this, we define the relation $\overline{E_L}$ as $\cl(E_L)$. $\overline{E_L}$ is $\emptyset$-invariant and contains $E_L$. There exists the finest bounded $\bigwedge$-definable over $\emptyset$ equivalence relation, denoted by $E_{KP}$ and known as equality of Kim-Pillay strong types (there is also an appropriate group of automorphisms $\autfKP(\C)$ such that $E_{KP}(a,b)$ if and only if for some $f\in\autfKP(\C)$, $a=f(b)$). The next theorem describes some relationship between $E_{KP}, \Theta$ and $E_L$.

\begin{theorem}\cite[Corollary 2.6]{pillay} \label{th:lascar}  $E_{KP} = \Theta \circ \overline{E_L}$
\end{theorem}

An attempt to understand the proof of this theorem was a starting point of this paper. In particular it was puzzling what properties of $E_L, E_{KP}$ and $\Theta$ are responsible for the relationship described in Theorem \ref{th:lascar}. It turnes out that the important point here is that both $E_L$ and $E_{KP}$ are orbit equivalence relations with respect to some groups of automorphisms of $\C$. We elaborate on this in Section 2. We generalize Theorem \ref{th:lascar} there and give a new proof of it based on Theorem \ref{th:ziegler}. Also in Section 2 we generalize some results about Lascar, Kim-Pillay and Shelah strong types.

Section 3 contains a model-theoretic analysis of a structure $N=(M,X,\cdot)$, where $M$ is a given stucture and $X$ is affine copy of some group $G$ definable in $M$. We describe the group of automorphisms of $N$ as a semi-direct product of $G$ and the group of automorphisms of $M$. In particular we reduce the question of $G$-compactness of $N$ to the question of $\bigwedge$-definability of a certain subgroup $G_L$ of $G$. This motivates us to look for examples of $G$, where $G_L$ is not $\bigwedge$-definable.

In Section 4 we verify that $G_L$ is $\bigwedge$-definable in several cases, e.g. when $M$ is small or simple or $o$-minimal and $G$ is definable compact.

In Section 5 we provide an example where a subgroup of $G$, similar in some sense to $G_L$, is not $\bigwedge$-definable, and also an example of a group $G$ that is not $G$-compact.

We assume that the reader is familiar with basic notions of model theory. 

The results in Sections 2, 3 and 4 are due to the first author, the proof of Lemma \ref{lem:el}(1) and the examples in Section 5 are due to the second author.

\section{Orbit equivalence relations}

In this section $G$ is always a subgroup of $\aut(\C)$. We can consider the orbit equivalence relation $E_G$ defined as follows: $E_G(a,b)$ if and only if there is some $f\in G$ with $a=f(b)$, where $a$ and $b$ are tuples of elements of $\C$ of length $\leq |T|$, such tuples are called small. In this paper we consider $E_G$ as an equivalence relation on the sets of small tuples of elements of various sorts of $\C$. 

The results of this section are concerned with various properties of relations of the form $E_G$. Our motivation is based on the observation that almost all important equivalence relations in model theory (e.g. $E_L$, $E_{KP}$ and $E_{Sh}$) are of this form.

Some statements from the next proposition are probably well known (see \cite{pillay,lascar,newelski,ziegler}).

\begin{lemma} \label{lem:basic}
\begin{enumerate}
\item[(i)] Let $M$ be an arbitrary small model, then \[ G\cdot \aut(\C/M) = \{f\in\aut(\C) : E_G(M,f(M)) \}.\]

\item[(ii)] The relation $E_G$ is $\emptyset$-invariant on every sort if and only if for every small $M \prec \C$ and every $F \in \aut(\C)$ \[G \subseteq G^F \cdot \aut(\C/M).\] 
In particular if $G$ contains $\bigcup_{F\in \aut(\C)} \aut(\C/F[M])$ for some small $M$, then $E_G$ is $\emptyset$-invariant if and only if $G \lhd \aut(\C)$.

\item[(iii)] If $G$ has bounded index in $\aut(\C)$, then $E_G$ is bounded and $E_L\subseteq E_G$. If $E_G$ is $\emptyset$-invariant bounded $G\lhd \aut(\C)$ and $G$ contains $\aut(\C/M)$ for some small $M$, then $G$ has bounded index in $\aut(\C)$.


\item[(iv)] Let $j\colon \aut(\C) \longrightarrow \galL(T)$ be the quotient map and assume that $\autfL(\C)\subseteq G$.
	\begin{enumerate}
	\item $j[G]$ is closed in $\galL(T)$ if and only if $E_G$ is $\bigwedge$-definable over any small model. If $G \lhd \aut(\C)$, then $\bigwedge$-definability is over $\emptyset$.
	\item $j[G]$ is open in $\galL(T)$ if and only if $G=\aut(\C/e)$ for some $e\in \acleq(\emptyset)$ (i.e. $e=\overline{m}/F$ for some $\emptyset$-definable finite equivalence relation $F$ on some $\C^n, n<\omega$).
	\end{enumerate}
\end{enumerate}
\end{lemma}
\begin{proof} (i) Easy.

(ii) Without loss of generality we may work with small models, because every tuple $a$ may be extended to small model $M$. Take an arbitrary small $M\prec \C$, $g \in G$ and $F \in \aut(\C)$. Then $E_G(M,g(M))$. Assume that $E_G$ is $\emptyset$-invariant. Then $E_G(F(M),F(g(M)))$ holds, so for some $g'\in G$, $F(g(M)) = g'(F(M))$. Thus $F^{-1} \circ g'^{-1} \circ F \circ g \in \aut(\C/M)$, so $g \in g'^F \circ \aut(\C/M) \subseteq G^F \circ \aut(\C/M) $. The other implication is similar.

For the second statement of (ii) assume that $G \subseteq G^F \cdot \aut(\C/M)$. Then conjugating by $F^{-1}$ we obtain \[ G^{F^{-1}} \subseteq G \cdot \aut(\C/F[M]) = G, \] for an appropriate small model $M$.

(iii) If $G$ has bounded index in $\aut(\C)$, then there is a normal subgroup $H \lhd \aut(\C)$ of bounded index, with $H\subseteq G$ (an intersection of boundedly many conjugates of $G$). Thus $E_H$ is bounded and invariant, so $E_L \subseteq E_H \subseteq E_G$.

For the second statement we use (i) to conclude that $G = \aut(\C / \ulcorner M/E_G\urcorner)$. $G$ has bounded index, because $M/E_G$ has boundedly many conjugates.

(iv) Note that $j^{-1}[j[G]] = G\cdot \autfL(\C) = G$, thus $\mu[G]=\nu^{-1}[j[G]]$ (because $j=\nu\circ\mu$). 

(a) $\Rightarrow$: Let $M$ be an arbitrary small model. If $j[G]$ is closed in $\galL(T)$, then $\mu[G]=\nu^{-1}[j[G]]=\{\tp(M'/M): \Phi(M',M)\}$ for some type $\Phi(x,y)$ over $\emptyset$. We have that \[E_G(a,b)\  \Longleftrightarrow\ \left( \exists f\in\aut(\C) \right) \left( a=f(b) \wedge \Phi(f(M),M) \right),\] and thus $E_G$ is $\bigwedge$-definable over $M$: \[E_G(a,b)\  \Longleftrightarrow\ (\exists z)(\tp(b,M)=\tp(a,z)\wedge \Phi(z,M)).\]

$\Leftarrow$: There is a type $\Phi(x,y)$ over $M$ such that \[ E_G(a,b)\  \Longleftrightarrow\ \Phi(a,b).\] Since $\mu[G]=\nu^{-1}[j[G]]$ it is enough to prove that $\mu[G]$ is closed in $S_M(M)$. This is clear, because: \[\mu[G] = \{\tp(g(M)/M) : g\in G \} = \{\tp(M'/M) : \Phi(M',M)\}.\] 

(b) $\Rightarrow$: First we deal with the case where $G \lhd \aut(\C)$. Since $\galL(T)$ is a compact topological group, $j[G]$ has finite index in $\galL(T)$, hence it is closed. By (iva) $E_G$ is $\emptyset$-$\bigwedge$-definable. Also $G$ has finite index in $\aut(\C)$. It follows that $E_G$ has finitely many classes on $\tp(M)(\C)$ (the set of realisations of type $\tp(M)$) and from (i) we have $G = \aut(\C/(M/E_G))$. Hence there are a finite $\emptyset$-definable equivalence relation $F$ and $\overline{m} \subset M$ such that $G = \aut(\C/(\overline{m}/F))$.

Now we deal with the general case, so $G<\aut(\C)$ need not be normal. However, still $G$ has finite index in $\aut(\C)$. Hence there is a normal subgroup $H \lhd \aut(\C)$ contained in $G$ and such that $j[H]$ is open (an intersection of finitely many conjugates of $G$). We may apply the first case to $H$. We get an $e\in \acleq(\emptyset)$ such that $H = \aut(\C/e)$. An element $e$ has finitely many conjugates, so $e'=\ulcorner\{g\cdot e : g\in G \}\urcorner \in \acleq(\emptyset)$. Now it is obvious that $G = \aut(\C/e')$.

$\Leftarrow$: The subset $\nu^{-1}[j[G]]=\mu[G]=\{\tp(f(M)/M) : F(\overline{m}, f(\overline{m})), f\in\aut(\C)\}$ of $S_M(M)$ is clopen.
\end{proof}

\begin{problem} Consider an equivalence relation $E$ on sorts of $\C$ which is $\emptyset$-invariant. Then we can build the following growing sequence of $\emptyset$-invariant relations:
\begin{enumerate}
\item[(i)] $E_0=E$,

\item[(ii)] $E_1=\cl(E)$ in $S_{k+k}(\emptyset)$,

\item[(iii)] for $1\leq \alpha\in Ord$ let 
	\begin{itemize}
	\item $E_{\alpha+1}=\cl(\textrm{transitive closure of }E_{\alpha})$,
	\item if $\alpha\in Lim$, then $E_{\alpha}= \bigcup_{\lambda<\alpha}E_{\lambda}$.
	\end{itemize}
\end{enumerate}
Take $E_{\infty}=\bigcup_{\alpha\in Ord}E_{\alpha}$. Then clearly $E_{\infty}$ is the finest type definable equivalence relation which extends $E$, so we may ask the question: what is the first $\alpha_E$ for which $E_{\alpha_E}=E_{\infty}$? If $E=E_G$, where $\autfL(\C)\subseteq G \lhd \aut(\C)$, then from the next Theorem \ref{th:moje1}(ii) we conclude that $\alpha_E\leq 2$.
\end{problem}

It can be proved that $\autfKP(\C)=j^{-1}[\cl(\id_{\galL(T)})]$. Recall that $E_{KP}=E_{\autfKP(\C)}$ is the finest bounded $\bigwedge$-definable over $\emptyset$ equivalence relation. The next Theorem \ref{th:moje1}(i) generalizes this remark and Theorem \ref{th:lascar} to an arbitrary group of automorphisms containing $\autfL(\C)$.

\begin{theorem}\label{th:moje1}
Let $\autfL(\C)\subseteq G<\aut(\C)$ and consider $\overline{\overline{G}}=j^{-1}[\cl(j[G])]$. Then
\begin{enumerate}
\item[(i)] On each sort of $\C$ the relation $E_{\overline{\overline{G}}}$ is the finest bounded $\bigwedge$-definable over any small model equivalence relation which extends $E_G$.

\item[(ii)] If additionally $G\lhd\aut(\C)$, then \[E_{\overline{\overline{G}}} = \Theta \circ \overline{E_G},\] where $\overline{E_G}$ is $\cl(E_G)$ in $S_{k+k}(\emptyset)$.
\end{enumerate}
\end{theorem}
\begin{proof} (i) Let $E$ be a $\bigwedge$-definable over $M$ equivalence relation and $E_G\subseteq E$. Take an arbitrary $f\in \overline{\overline{G}}$ and a small tuple $b$. We have to prove that $E(f(b),b)$. Consider the following set \[H=\{f\in\aut(\C) : E(f(b),b)\}\] ($H$ is not necessarily a group, because $E$ is not necessarily $\emptyset$-invariant). It is enough to show that $\overline{\overline{G}}\subseteq H$. 

Note that $j^{-1}[j[H]]=\autfL(\C)\cdot H=H$, because for $f\in\autfL(\C)$, $h\in H$ we have $E(h(b),b)$ and $E(f(h(b)),h(b))$ ($E_L\subseteq E$), so $E(f(h(b)),b)$ and $f\circ h\in H$. 

Since $E_G\subseteq E$ we have $G\subseteq H$, so we must only prove that $\cl(j[G])\subseteq j[H]$ (because $j^{-1}[j[H]] = H$). The proof is completed by showing that $j[H]$ is closed in $\galL(T)$. This follows from the fact that the set \[\nu^{-1}[j[H]]=\mu[H]=\{\tp(f(M')/M') : E(f(b),b), f\in\aut(\C)\}\] is closed in $S_{M'}(M')$, where $Mb\subseteq M' \prec \C$.

(ii) The relation $E_{\overline{\overline{G}}}$ is $\bigwedge$-definable over $\emptyset$, so $E_{\overline{\overline{G}}}$ is a closed subset of $S_{k+k}(\emptyset)$, thus $\overline{E_G}\subseteq E_{\overline{\overline{G}}}$. This gives $\Theta \circ \overline{E_G}\subseteq E_{\overline{\overline{G}}}$.

Now we prove that $E_{\overline{\overline{G}}} \subseteq \Theta \circ \overline{E_G}$. Assume that $a,b$ are small tuples such that $E_{\overline{\overline{G}}}(a,b)$, i.e. $a=f(b)$ for some $f\in \overline{\overline{G}}$. Without loss of generality we may assume that $b=M$, for some small $M\prec\C$, so $a=f(M)$. Let $p=\mu(f)=\tp(f(M)/M)$. Then $\nu(p)=j(f)\in\cl(j[G])$ and \[[p]_{\Theta} \cap \cl(\nu^{-1}[j[G]])\neq\emptyset,\] because otherwise $[p]_{\Theta}\subseteq \inter(\nu^{-1}[j[G]^{c}])$, and from Theorem \ref{th:ziegler} \[\nu(p)\in\inter(\nu[\nu^{-1}[j[G]^{c}]]) =\inter(j[G]^{c}) = \cl(j[G])^{c},\] a contradiction. 

Let $q\in [p]_{\Theta} \cap \cl(\nu^{-1}[j[G]])$. There is some $c=M'\models q$ such that $\Theta(f(M),M')$ and $q=\tp(M'/M)$ is in 
\begin{equation}
\begin{array}{rl}
\cl(\nu^{-1}[j[G]]) = \cl(\mu[G]) &= \cl\{\tp(g(M)/M) : g\in G \} \nonumber \\
&=\cl\{\tp(g(M)/M) : E_G(g(M),M)\}. \nonumber
\end{array}
\end{equation}
Finally $\tp(M',M)\in\cl(E_G)=\overline{E_G}$, and we obtain that $\Theta(a,c)$ and $\overline{E_G}(c,b)$.
\end{proof}

Now we consider the relation $E_{Sh}$ of equality of Shelah strong types: \[E_{Sh} = \bigcap\{E : E\textrm{ is a }\emptyset\textrm{-definable finite equivalence relation}\}.\] It can be proved that $E_{Sh} = E_{j^{-1}[\QC]}$, where $\QC \lhd \galL(T)$ is the intersection of all open subgroups of $\galL(T)$ (the quasi-connected component). When $\galL(T)$ is Hausdorff (i.e. $T$ is $G$-compact) then $\QC$ is just the connected component of $\galL(T)$.

In the next proposition we generalize this property of $E_{Sh}$, but first we need a definition: if $A\subseteq \galL(T)$, then by $\QC(A)$ we denote the following set \[\bigcap\{H<\galL(T) : A\subseteq H\textrm{ and }H\textrm{ is open}\}.\]  

\begin{proposition} \label{prop:sh}
If $H< \galL(T)$, then $E_{j^{-1}[\QC(H)]}$ is the intersection of all $\emptyset$-definable finite equivalence relations which extend $E_{j^{-1}[H]}$: \[E_{j^{-1}[\QC(H)]} = \bigcap\{E : E\textrm{ is a }\emptyset\textrm{-definable finite e.r. and } E_{j^{-1}[H]}\subseteq E\}.\] Moreover $j^{-1}[\QC(H)]$ is equal to the group of all $f\in\aut(\C)$, satisfying \[E_{j^{-1}[\QC(H)]}(a,f(a))\] for arbitrary small tuple $a$.
\end{proposition}
\begin{proof} First we prove the equality of relations. ($\subseteq$) Assume that small tuples $a,b$ are $E_{j^{-1}[\QC(H)]}$ equivalent, so $a=f(b)$ for some $f\in j^{-1}[\QC(H)]$, and $E$ is a $\emptyset$-definable finite equivalence relation extending $E_{j^{-1}[H]}$. Define \[G' = \{f\in\aut(\C) : E(f(b),b)\}=\aut(\C/(b/E)).\] Then $H\subseteq j[G']$ and $j[G']$ is open as a subset of $\galL(T)$ (Lemma \ref{lem:basic}(iv)(b)). Therefore $\QC(H)\subseteq j[G']$ and \[f\in j^{-1}[\QC(H)] \subseteq j^{-1}[j[G']]=G'\cdot\autfL(\C)=G',\] so $E(a,b)$ holds.

($\supseteq$) Let $\QC(H)=\bigcap\{G_i : i\in I\}$. Using Lemma \ref{lem:basic}(iv)(b) we can find $(e_i)_{i\in I}\subseteq \acleq(\emptyset)$ such that $j^{-1}[G_i]=\aut(\C/e_i)$. Then $j^{-1}[QC(H)]=\aut(\C/\{e_i\}_{i\in I})$. We can assume that $e_i=m_i/F_i$ for some $\emptyset$-definable finite equivalence relations $F_i$. Assume that $(a,b)$ belongs to \[ \bigcap\{E : E\textrm{ is a }\emptyset\textrm{-definable finite e.r. and } E_{j^{-1}[H]}\subseteq E\}.\] We have to find $f\in j^{-1}[QC(H)]$ for which $b=f(a)$. It suffices to prove that the following type in variables $(y_i)_{i\in I}$ is consistent:
\[ ``\tp(b,y_i)_{i\in I} = \tp(a,m_i)_{i\in I}" \wedge \bigwedge_{i\in I} F_i(m_i,y_i).\] Let $\varphi(x,x_1,\ldots,x_n)\in \tp(a,m_i)_{i\in I}$. It is enough to show that \[\psi(b,m_1,\ldots,m_n)=(\exists y_1,\ldots,y_n) \left(\varphi(b,y_1,\ldots,y_n)\wedge \bigwedge_{1\leq i \leq n}F_i(m_i,y_i)\right)\] holds. The formula $\psi=\psi(x,m_1,\ldots,m_n)$ is almost over $\emptyset$ (because $m_1,\ldots,m_n\in \acl(\emptyset)$). 
Let $\psi_1,\ldots,\psi_k$ be all conjugates of $\psi$ over $\emptyset$ and take \[A(x,y)=\bigwedge_{1\leq i \leq k} (\psi_i(x)\leftrightarrow\psi_i(y)).\] $A$ is a $\emptyset$-definable finite equivalence relation and $E_{j^{-1}[H]}\subseteq A$\  (because $j^{-1}[H] \subseteq j^{-1}[QC(H)] = \aut(\C/(m_i/F_i))$). Therefore $A(a,b)$ and we know that $\psi(a,m_1,\ldots,m_n)$ holds, so $\psi(b,m_{\leq n})$ also holds.

Now we prove the second part of the proposition. Let $G'$ be the group of all automorphisms preserving $E_{j^{-1}[\QC(H)]}$. Inclusion $j^{-1}[\QC(H)] \subseteq G'$ is obvious. 

($\supseteq$) Let $g\in G'$ and $a=M$ be a small model. Then $E_{j^{-1}[\QC(H)]}(M,g(M))$, so $g(M) = f(M)$ for some $f\in j^{-1}[\QC(H)]$. Thus $j(g) = j(f)$ (because $gf^{-1}\in \autfL(\C)$) and $j(f)\in \QC(H)$. Therefore $g\in j^{-1}[\QC(H)]$.
\end{proof}



\section{An Example}

Let $M$ be an arbitrary structure in which we have a $\emptyset$-definable (interpretable) group $G$. In this section we consider the following two sorted structure: $N=(M,X,\cdot)$, where
\begin{itemize}
\item $X$ and $M$ are disjoint sorts,
\item $\cdot\colon G\times X \rightarrow X$  is a regular (free and transitive) action of $G$ on $X$ i.e. $X$ is an affine copy of $G$,
\item on $M$ we take its original structure.
\end{itemize}

This structure was already considered e.g. in \cite{ziegler} and \cite{newelski}. Our study of $N$ is based on ideas from \cite[Section 7]{ziegler}. 

In this section we describe various groups of automorphisms of $N$ in terms of appropriate groups of automorphisms of $M$ and groups related to $G$. We also give a description of the relations $E_L$, $E_{KP}$ and $E_{Sh}$ on the sort $X$ of $N$. In particular, in Corollary \ref{cor:comp} we prove that $G$-compactness of $N$ is equivalent to $G$-compactness of $M$ and $\bigwedge$-definability of a certain subgroup $G_L$ of $G$. Thus constructing a group $G$ where the subgroup $G_L$ is not $\bigwedge$-definable may yield a new example of a non-$G$-compact theory.

Fix an arbitrary point $x_0$ from $X$ and take $N^*=(M^*,X^*,\cdot)$, a monster model extending $N$. Then $G\subseteq G^*$ and $X=G\cdot x_0 \subseteq G^*\cdot x_0=X^*$.

The group $G^*$ acts on itself in two different, but commuting ways, the first one is by left translation $(g,h) \mapsto gh$, and the second one by the following rule $(g,h) \mapsto hg^{-1}$. We define homomorphic embeddings of automorphism groups: 
\[ \overline{\cdot}\colon \aut(M^*) \hookrightarrow \aut(N^*),\ \  \overline{\cdot}\colon G^* \hookrightarrow \aut(N^*).\] Let $h\in G^*, f\in \aut(M^*), g\in G^*$. We define $\overline{f}, \overline{g} \in \aut(N^*)$ by: \[ \overline{f}|_{M^*}=f,\  \overline{f}(h\cdot x_0)=f(h)\cdot x_0,\]
\[ \overline{g}|_{M^*}=\id_{M^*},\  \overline{g}(h\cdot x_0)=(hg^{-1})\cdot x_0.\] 
It is easy to verify the following laws: for $f\in \aut(M^*)$, $g\in G^*$ we have  \[\overline{f}\circ \overline{g} = \overline{f(g)}\circ \overline{f},\ \ \overline{g}\circ \overline{f} = \overline{f}\circ \overline{f^{-1}(g)}.\]
Using these embeddings we can identify $\aut(M^*)$ and $G^*$ with their images in $\aut(N^*)$ and conclude that $G^*\lhd \aut(N^*)$. In fact we will prove that $\aut(N^*)$ is a semi-direct product of $G^*$ and $\aut(M^*)$. 

There are two different actions of the group $G^*$ on the set $X^*$: the first one comes from the above embedding \[\overline{g}(h\cdot x_0) = (hg^{-1})\cdot x_0\] (it is definable over $x_0$). The second one comes from the regular action \[g\cdot (h\cdot x_0) = (gh)\cdot x_0.\] If $A\subseteq G^*$ satisfies $hA^{-1} = Ah$, then the orbits of $h\cdot x_0$ under  both actions coincide: $\overline{A}(h\cdot x_0) = A\cdot(h\cdot x_0)$ (in this case we just write $A\cdot(h\cdot x_0)$).

In order to describe properties of $N^*$ in terms of $M^*$ and $G^*$ we need the next definition.

\begin{definition} \label{def:kom}
For a group $G$ and a binary relation $E$ on $G$ we define the set of $E$-commutators $X_E = \{a^{-1}b : a,b\in G, E(a,b)\}$ and the $E$-commutant $G_E$ as the subgroup of $G$ generated by $X_E$ \[G_E=\langle X_E \rangle <G.\]
\end{definition}

\begin{remark}
If $E=E_H$ for some $H<\aut(G,\cdot)$, then $G_{E_H}\lhd G$. If $E$ is $\emptyset$-invariant, then $X_E$ and $G_E$ are also $\emptyset$-invariant. If $E$ is bounded, then $G_E$ has bounded index in $G$, moreover $[G : G_E]\leq |G/E|$.
\end{remark}
\begin{proof}
Let $a,x \in G$ and $h\in H$. Then \[(X_{E_H})^x \ni(a^{-1}h(a))^x =(ax)^{-1}h(a)x = ((ax)^{-1}h(ax)) (h(x)^{-1}x)\in X_{E_H}^2.\]  The last statement follows from the observation: if $a^{-1}b\notin G_E$, then $\neg E(a,b)$.
\end{proof}

The following example justifies the names ``$E$-commutators" and ``$E$-commutant" from the previous definition. Let $E$ be the conjugation relation in $G$ i.e. $E=E_{\inn(G)}$ (where $\inn(G)$ is the group of inner automorphisms of $G$). Then $X_E$ is the set of all commutators and $G_E=[G,G]$. 

In the case where $E=E_L$ [$E=E_{KP}, E_{Sh}$, respectively] we just write $X_L$ and $G_L$ [$X_{KP}$, $X_{Sh}$] instead of $X_{E_L}$ and $G_{E_L}$ [$X_{E_{KP}}$, $X_{E_{Sh}}$]. Note that $G_L$ is generated by $X_{\Theta}$. 

In the next proposition we describe $\aut(N^*), \autfL(N^*)$ and $\galL(\Th(N))$ as semidirect products of automorphisms groups of $M^*$ and appropriate groups associated with $G$.

\begin{proposition} \label{prop:1}
\begin{enumerate}
\item $\aut(N^*) = G^* \rtimes \aut(M^*)$, more precisely: for $F\in\aut(N^*)$, $F=\overline{g}\circ\overline{f}$, where $f=F|_{M^*}$ and $F(x_0)=g^{-1}\cdot x_0$.

\item Let $(N',X')\prec(N^*,X^*)$ and $X'=G'\cdot(h_0\cdot x_0) $ for some $h_0\cdot x_0 \in X'$. Then \[F\in\aut(N^*/N')\  \Longleftrightarrow\  (\exists f\in\aut(M^*/M')) \left(F= \overline{f} ^ {\overline{h_0}}\right).\]

\item $\autfL(N^*) = G^*_L \rtimes \autfL(M^*)$ and $\galL(\Th(N)) = G^*/G^*_L \rtimes \galL(\Th(M))$.
\end{enumerate}
\end{proposition}

\begin{proof} (1) Let $F\in\aut(N^*)$ and $f=F|_{M^*}$. Then $F\overline{f}^{-1}$ is the identity on $M^*$, and on $X^*=G^*\cdot x_0$ we have: \[F\overline{f}^{-1}(h\cdot x_0) = F(f^{-1}(h)\cdot x_0) = h\cdot F(x_0) = h\cdot(g^{-1}\cdot x_0) = \overline{g}(h\cdot x_0),\] for some $g\in G^*$. Thus $F=\overline{g}\circ\overline{f}$. The group $\aut(M^*)$ acts on $G^*$ by conjugation, so for $g\in G^*$ and $f\in\aut(M^*)$, $\overline{g}^{\overline{f}} = \overline{f(g)} \in G^*$. It is clear that $G^*\cap\aut(M^*) = \{0\}$.

(2) ($\Leftarrow$)  It is clear that $F|_{M'}=\id_{M'}$. Using the fact that $f|_{M'} = \id_{M'}$ we get for $h'\in X'$: \[ \overline{f}^{\overline{h_0}}(h'h_0\cdot x_0)=\overline{h_0^{-1}}\circ \overline{f}(h'\cdot x_0) = \overline{h_0^{-1}}(h'\cdot x_0) = h'h_0\cdot x_0.\] Thus $\overline{f}^{\overline{h_0}}|_{X'}=\id_{X'}$.

($\Rightarrow$) Let $f = F|_{M^*}$. Then $f = \overline{f}^{\overline{h_0}}|_{M^*} \in \aut(M^*/M')$. By assumptions \[ h_0\cdot x_0=F(h_0\cdot x_0) = F(h_0)\cdot F(x_0) = f(h_0)\cdot F(x_0),\] and then $F(x_0) = f(h_0^{-1})h_0\cdot x_0$. By (1), $F=\overline{h_0^{-1} f(h_0)}\circ\overline{f} = \overline{h_0^{-1}}\circ\overline{f}\circ\overline{h_0} = \overline{f}^{\overline{h_0}}$.

(3) It suffices to prove the first equality. $\subseteq$: From (2) we conclude that for every $F\in\autfL(N^*)$ there are $h_1,\dots,h_n\in G^*$ and $f_1,\dots,f_n\in\autfL(M^*)$ such that $F = \overline{f_1}^{\overline{h_1}}\circ\ldots\circ\overline{f_n}^{\overline{h_n}}$. Then \[F = \overline{h_1^{-1}f_1(h_1)}\circ\overline{f_1} \circ \overline{h_2^{-1}f_2(h_2)}\circ\overline{f_2} \circ \ldots \circ \overline{h_n^{-1}f_n(h_n)}\circ\overline{f_n}.\] Using the rule $\overline{f}\circ \overline{g} = \overline{f(g)}\circ \overline{f}$, one can prove that $F = \overline{g}\circ\overline{f_1\ldots f_n}$, for some $g\in G_L$ (for example 
$\overline{f_1} \circ \overline{h_2^{-1}f_2(h_2)} =  \overline{f_1(h_2^{-1})f_1(f_2(h_2))} \circ \overline{f_1}$, and $f_1(h_2^{-1})f_1(f_2(h_2)) = f_1(h_2)^{-1} f_2^{f_1^{-1}} (f_1(h_2)) \in X_L$).

$\supseteq$: It is clear that $\autfL(N^*) \supseteq \autfL(M^*)$ (use (2)). It is enough to prove that $\autfL(N^*) \supseteq X_L$. Assume that small tuples $a,b$ satisfy $b=f(a)$, for some $f\in\autfL(M^*)$. We have to prove that $\overline{a^{-1}b}\in\autfL(N^*)$. Since $\overline{f}^ {\overline{a}}\in \autfL(N^*)$, we have $\overline{a^{-1}b} = \overline{a^{-1}f(a)} = \overline{f}^ {\overline{a}}\circ\overline{f^{-1}} \in\autfL(N^*)$.
\end{proof}

Now we characterize some invariant subgroups of $G^*$: $G^{0}_{\emptyset}, G^{00}_{\emptyset}$ and $G^{\infty}_{\emptyset}$, in terms of $N^*$.

\begin{proposition} \label{prop:1:2}
\begin{enumerate}
\item $G^*_L = G^*\cap \autfL(N^*)$ and $G^*_L$ is the smallest $\emptyset$-invariant subgroup of $G$ with bounded index in $G^*$ (i.e. $G^*_L = G^{\infty}_{\emptyset}$).

\item Let $G'_{KP} = G^*\cap\autfKP(N^*)$, then $G^*_{KP} \subseteq G'_{KP}$ and $G'_{KP}$ is the smallest $\bigwedge$-definable over $\emptyset$ subgroup with bounded index in $G^*$ (i.e. $G'_{KP} = G^{00}_{\emptyset}$).

\item Let $G'_{Sh} = G^*\cap\autfSh(N^*)$, then $G^*_{Sh}\subseteq G'_{Sh}$ and $G'_{Sh}$ is the intersection of all $\emptyset$-definable subgroups of $G^*$ with finite index (i.e. $G'_{Sh} = G^{0}_{\emptyset}$).
\end{enumerate}
\end{proposition}
\begin{proof}
(1) The first equality follows directly from Proposition \ref{prop:1}(3). Let $H<G^*$ be $\emptyset$-invariant with bounded index. It suffices to prove that $X_{\Theta}\subseteq H$. Take an order inscernible sequence $(a_n)_{n<\omega}$ (so $\Theta(a_0,a_1)$). If $a_0^{-1}a_1\notin H$, then for every $i<j<\omega,\ a_i^{-1}a_j \notin H$, but we can extend an indiscernible sequence as much as we want, so the index $[G^* : H]$ is unboundedly large, a contradiction.

(2) If $N'\prec N^*$ is an arbitrary small model, then \[G'_{KP} = \{g\in G^* : E_{KP}(N', \overline{g}(N'))\}.\] Inclusion $\subseteq$ is obvious. $\supseteq$: If $E_{KP}(N', \overline{g}(N'))$, then $\overline{g}(N') = F(N')$ for some $F\in \autfKP(N^*)$. Since $\overline{g}|_{M^*} = \id_{M^*}$, $F\in G^*$, and $\overline{g} F^{-1} \in \autfL(N^*)$, so $\overline{g} \in G'_{KP}$.

$G'_{KP}$ has bounded index ($G^*_L\subseteq G'_{KP}$) and is $\bigwedge$-definable over $N'x_0$. In fact $G'_{KP}$ is $\emptyset$-invariant. To see this let $F = \overline{g'}\circ\overline{f'} \in\aut(N^*)$. Then 
\[F[G'_{KP}] = \{ F(g) : E_{KP}(N', \overline{g}(N')) \} = \{ F(g) : E_{KP}(F(N'), F(\overline{g}(N'))) \}, \] but $F\circ\overline{g} = \overline{g'}\circ\overline{f'}\circ\overline{g} = \overline{g'f'(g)g'^{-1}} \circ F$, thus 
\begin{equation}
\begin{array}{rl}
F[G'_{KP}] &= \{ f'(g) : E_{KP}(F(N'), \overline{g'} \circ \overline{f'(g)g'^{-1}} (F(N'))) \} \nonumber \\
&= \{ f'(g) : E_{KP}(\overline{g'^{-1}}F(N'), \overline{f'(g)}(\overline{g'^{-1}}F(N')) \}  = G'_{KP}, \nonumber
\end{array}
\end{equation}
and hence $G'_{KP}$ is $\bigwedge$-definable over $\emptyset$. 
The relation $E(x,y)=x^{-1}y\in G'_{KP}$ is bounded $\bigwedge$-definable over $\emptyset$, therefore $E_{KP}|_{G^*} \subseteq E$ and $G^*_{KP}\subseteq G'_{KP}$. Take $H<G^*$, another subgroup which is $\bigwedge$-definable over $\emptyset$ and has bounded index in $G^*$. Then $E_{KP}\subseteq E_H$, so for $g\in G'_{KP}$ we have $E_H(x_0,\overline{g}(x_0))$ and then $g^{-1}\cdot x_0=\overline{g}(x_0) = \overline{h}(x_0) = h^{-1}\cdot x_0$ for some $h\in H$. By regularity of $\cdot$ we obtain $g=h \in H$.

(3) As in (2) it can be proved that $G'_{Sh}$ is $\bigwedge$-definable over $\emptyset$. 
Let $g\in G'_{Sh}$, and $H< G^*$ be a $\emptyset$-definable subgroup with finite index in $G^*$. We show that $g\in H$. Consider the relation $E(x,y)= (\exists h\in H)(x = h\cdot y)$ on $X^*$. $E$ is a $\emptyset$-invariant, finite equivalence relation on $X^*$, thus $E_{Sh}|_{X^*} \subseteq E$. By regularity of $\cdot$ we conclude that $g\in H$. If we consider $E(x,y) = x^{-1}y\in H$ on $G^*$, then $E_{Sh}|_{G^*} \subseteq E$ and therefore $G^*_{Sh} \subseteq H$.

Let $g$ belong to all $\emptyset$-definable subgroups of $G^*$ of finite index. We prove that $\overline{g}\in \autfSh(N^*)$. From Proposition \ref{prop:sh} we know that $\autfSh(N^*)$ is the preimage under the quotient map $j$ of the quasi-connected component $\QC$ of $\galL(\Th(N))$. Let $H \lhd \galL(\Th(N))$ be an open subgroup. It suffices to show that $\overline{g} \in j^{-1}[H] \lhd \aut(N^*)$. Note that the group $H' = j^{-1}[H]\cap G^*$ is $\emptyset$-invariant, because for $f\in\aut(M^*)$ if $\overline{g}\in j^{-1}[H]$, then $\overline{f(g)} = \overline{g}^{\overline{f^{-1}}} \in j^{-1}[H]$. $H'$ is also definable, because by Lemma \ref{lem:basic}(iv)(b), $j^{-1}[H] = \aut(\C/\overline{m}/F))$, so $g\in H'$ if and only if $F(\overline{m},\overline{g}(\overline{m}))$). Hence $H'$ is a $\emptyset$-definable subgroup of $G^*$ of finite index and thus $\overline{g} \in j^{-1}[H]$.
\end{proof}

The compact topological group $\galL(\Th(N^*))$ contains as a subgroup the group $G^*/G^*_L$, so we may ask about the induced topology on $G^*/G^*_L$. The next proposition describes this topology.

\begin{proposition} \label{prop:2}
\begin{enumerate}
\item The induced subspace topology on $G^*/G^*_L$ from $\galL(\Th(N))$ is precisely the logic topology: let $i\colon G^* \rightarrow G^*/G^*_L$ be the quotient map, then $X\subseteq G^*/G^*_L$ is closed if and only if its preimage $i^{-1}[X] \subseteq G^*$ is $\bigwedge$-definable over some (equivalently every) small model. With this topology $G^*/G^*_L$ is a compact topological group (this topology is Hausdorff if and only if $G^*_L$ is $\bigwedge$-definable).

\item The topology of $\galL(\Th(M))$ as the Lascar group of $\Th(M)$ and the induced topology on $\galL(\Th(M))$ as a subspace of $\galL(\Th(N))$ coincide.

\item If $X\subseteq G^*/G^*_L$ and $Y\subseteq \galL(\Th(M))$ are closed, then $X\cdot Y\subseteq \galL(\Th(N))$ is also closed. In particular, if $\Th(M)$ is $G$-compact, then $G^*/G^*_L$ is closed subgroup of $\galL(\Th(N))$.

\item The closure of identity in $G^*/G^*_L$ is $G'_{KP}/G^*_L$.

\item The quasi-connected component (the intersection of all open subgroups) of $G^*/G^*_L$ is $G'_{Sh}/G^*_L$.
\end{enumerate}
\end{proposition}
\begin{proof}
(1) Let $N'$ be a small model. Without loss of generality we may assume that $x_0 \in N'$. The restriction of the quotient map $j$ to $G^*$ is precisely the quotient map $i$. We have the following commutative diagram:
\[\xymatrix{
G^* \ar@{>>}[r]^-{i} \ar@{^{(}->}[d]^-{\subseteq} & G^*/G^*_L \ar@{^{(}->}[d]_-{\subseteq} \\
\aut(N^*) \ar@{>>}[r]^-{j} & \galL(\Th(N))
}\]
Let $X\subseteq G^*/G^*_L$ be closed in the induced subspace topology, i.e. $X = G^*/G^*_L \cap C$, where $C \subseteq \galL(\Th(N))$ is closed. Then $\nu^{-1}[C]$ is closed in $S_{N'}(N')$, so there exists a type $\Phi(x,y)$ over $\emptyset$ for which
\[ \nu^{-1}[C] = \mu[j^{-1}[C]] = \{ \tp(F(N')/N') : F \in j^{-1}[C] \} = 
\{ \tp(N''/N') : \Phi(N'',N') \}.\] The subset $i^{-1}[X] \subseteq G^*$ is $\bigwedge$-definable over $N'$, because for $g \in G^*$
\[ g\in i^{-1}[X]\  \Leftrightarrow\  \overline{g} \in j^{-1}[C]\  \Leftrightarrow\  \Phi(\overline{g}(N'),N').\] 
The implication $\Leftarrow$ in the last equivalence holds, because if $\Phi(\overline{g}(N'),N')$, then $\overline{g}(N') = F(N')$ for some $F \in j^{-1}[C]$, and thus $j(\overline{g}) = j(F) \in C$.

Now assume that $i^{-1}[X]$ is $\bigwedge$-definable over $N'$, i.e. for $g\in G^*,\ g\in i^{-1}[X]$ if and only if $\Psi(g,N')$, for some type $\Psi$. Let $C = X \cdot \galL(\Th(M)) \subseteq \galL(\Th(N))$. Then \[X = G^*/G^*_L \cap C.\] In order to prove that $C$ is closed in $\galL(\Th(N))$ it is enough to show that 
\begin{equation}
\begin{array}{rl}
\nu^{-1}[C] &= \{ \tp(F(N')/N') : F\in j^{-1}[C] \} \nonumber \\
 &= \{ \tp(N''/N') : x_0^{N''} = g^{-1}\cdot x_0,\ \Psi(g,N') \textrm{ holds and } \tp(N'') = \tp(N')\}. \nonumber
\end{array}
\end{equation}
The last equality holds because $j^{-1}[C] = i^{-1}[X] \circ \aut(M^*)$, and if $F = \overline{g}\circ \overline{f},\ g\in i^{-1}[X]$, then $x_0^{N''} = F(x_0) = \overline{g}\circ \overline{f}(x_0) = g^{-1}\cdot x_0$ (here $N'' = F(N')$).

(2) The proof is similar to the proof in (1) and we leave it to the reader.

(3) The set $\nu^{-1}[X\cdot Y] = \mu[j^{-1}[X\cdot Y]]$ is closed in $S_{N'}(N')$ because it is equal to the following 
\begin{equation}
\begin{array}{ll}
\{ \tp(\overline{g} \circ \overline{f}(N')/N') : g\in i^{-1}[X], f\in j^{-1}[Y] \} = \nonumber \\
\{ \tp(N''/N') : \tp(M''/M')\in \nu^{-1}[Y],\ x_0^{N''} = g^{-1}\cdot x_0,\ g\in i^{-1}[X]\textrm{ and }  \tp(N'') = \tp(N')\}. \nonumber
\end{array}
\end{equation}
Above we use the fact that $j^{-1}[X\cdot Y] = i^{-1}[X] \circ \nu^{-1}[Y]$.

(4) $G'_{KP}/G^*_L$ contains $\cl(\id)$, because $G'_{KP}$ is $\bigwedge$-definable over $\emptyset$. The subgroup $i^{-1}[\cl(\id)]$ of $G^*$ is $\bigwedge$-definable over $\emptyset$ and of bounded index (because $G^*_L\subseteq i^{-1}[\cl(\id)]$), thus $G'_{KP}\subseteq i^{-1}[\cl(\id)]$.

(5) The group $G'_{Sh}$ is the intersection of all $\emptyset$-definable subgroups of $G^*$ of finite index, thus $G'_{Sh}/G^*_L$ contains quasi-connected component $\QC$ (because if $H<G^*$ is $\emptyset$-definable of finite index, then $H/G^*_L$ is closed of finite index, hence open). Let $H$ be an arbitrary open subgroup of $G^*/G^*_L$. It suffices to show that $G'_{Sh}/G^*_L \subseteq H$. The group $H$ is closed of finite index, hence $H\cdot \galL(\Th(M))$ is a closed subgroup of $\galL(\Th(N))$ of finite index. Therefore \[\autfSh(N^*) \subseteq j^{-1}[H\cdot \galL(\Th(M))],\] and then $G'_{Sh} \subseteq i^{-1}[H]$. This gives $G'_{Sh}/G^*_L \subseteq H$.
\end{proof}

The next corollary motivates us to investigate $\bigwedge$-definability of $G^*_L$. We do this in the next section. If $G^*_L$ is not $\bigwedge$-definable, then $N$ may give us a new kind of not $G$-compact theory.

\begin{corollary} \label{cor:comp}
$\Th(N)$ is $G$-compact if and only if $\Th(M)$ is $G$-compact and $G^*_L$ is $\bigwedge$-definable.
\end{corollary}
\begin{proof} The topological group $G$ is Hausdorff if and only if $\{e_G\}$ is closed and we can apply the previous proposition.
\end{proof}

Now we describe the relations $\Theta, E_L, E_{KP}$ and $E_{Sh}$ on the sort $X^*$ in terms of orbits of the groups $G^*_L, G'_{KP}$ and $G'_{Sh}$ from Proposition \ref{prop:1}.

\begin{lemma} Let $x\in X^*$ and $n<\omega$. \label{lem:el}
\begin{enumerate}
\item $\{y\in X^* : \Theta^n(x,y)\}=X_{\Theta}^n\cdot x$ 
\item $x/E_L = G^*_L\cdot x$
\item $x/E_{KP}=G'_{KP} \cdot x$
\item $x/E_{Sh}=G'_{Sh} \cdot x$
\end{enumerate}
\end{lemma}
\begin{proof}
(1) It is enough to prove this for $n=1$. $\subseteq$: Assume $x,y\in X^*, \Theta(x,y)$ and $y=g_0x$ for some $g_0\in G^*$. We may assume that $x=x_0$. We can extend $(x_0,g_0x_0)$ to an order indiscernible sequence $(x_0, g_0x_0, g_1x_0,\ldots )\subseteq X^*$. Then for $0\leq i_1<\ldots<i_n<\omega, 0\leq j_1<\ldots<j_n<\omega$: 
\[(x_0, g_{i_1}x_0, g_{i_2}x_0, \ldots)\equiv (g_{j_1}x_0, g_{j_2}x_0, g_{j_3}x_0, \ldots).\] Applying the automorphism $\overline{g_{j_1}}$ we obtain: \[(g_{j_1}x_0, g_{j_2}x_0,g_{j_3}x_0,\ldots)\equiv (x_0, g_{j_2}g_{j_1}^{-1}x_0, g_{j_3}g_{j_1}^{-1}x_0,\ldots).\] Hence from the previous two equivalences we get
\[(g_{i_1}x_0, g_{i_2}x_0,\ldots)\underset{x_0}{\equiv} (g_{j_2}g_{j_1}^{-1}x_0, g_{j_3}g_{j_1}^{-1}x_0,\ldots),\] so
\[(g_{i_1}, g_{i_2},\ldots)\equiv (g_{j_2}g_{j_1}^{-1}, g_{j_3}g_{j_1}^{-1},\ldots).\] It means that $(g_0,g_1,\ldots)\subseteq G^*$ is also order indiscernible and $g_0\equiv g_0g_1^{-1}$, so $g_0 \in X_{\Theta}$.

$\supseteq$: Let $y=gx_0$ for $g=ab^{-1}\in X_{\Theta}$, where $\Theta(a,b)$. We can find an indiscernible sequence $(b,gb,\ldots)\subseteq G^*$, and then $(bx_0,gbx_0,\ldots)\subseteq X^*$ is also indiscernible, so $\Theta(bx_0,gbx_0)$. Applying $\overline{b}$, we obtain $\Theta(x_0,gx_0)$. 

(2) Inclusion $\supseteq$ follows from Proposition \ref{prop:1}(3). $\subseteq$: Let $y=F(x)$ for some $F=\overline{g}\circ\overline{f}\in\autfL(N^*)$. We may assume that $x=x_0$. Then $y=\overline{gf}(x_0)=\overline{g}{x_0}=g^{-1}x_0$ and $g\in G_L$.

(3) $\supseteq$ follows from Proposition \ref{prop:1:2}(2). Since $E_{KP}|_{X^*}\subseteq E_{G'_{KP}}|_{X^*}$ we have $\subseteq$.

(4) $\supseteq$ follows from Proposition \ref{prop:1:2}(3). $\subseteq$: We know that $G'_{Sh} = \bigcap_{i\in I} H_i$, where $H_i$ is $\emptyset$-definable with finite index. Therefore $E_{G'_{Sh}}|_{X^*} = \bigcap_{i\in I} E_{H_i}|_{X^*}$, so $E_{Sh}|_{X^*}\subseteq E_{G'_{Sh}}|_{X^*}$ and we are done.
\end{proof}

Using Theorems 1.1 and 3.1 from \cite{newelski} we can give a detailed analysis of Lascar and Kim-Pillay strong types on $X^*$. This analysis describes also some basic properties of the group $G$. By $\diam(a)$ we denote the diameter of the Lascar strong type $a/E_L$ (see \cite{newelski}). Note that every two elements of $X^*$ have the same type over $\emptyset$, thus their Lascar strong types have the same diameter. 

\begin{remark} There are only two possibilities: \label{rem:case}
\begin{enumerate}
\item The diameters of all Lascar strong types on $X^*$ are infinite. The group $G^*_L$ is not $\bigwedge$-definable, $E_L \subsetneqq E_{KP}$, $G^*_L \subsetneqq G'_{KP}$ (i.e. $\Th(N)$ is not $G$-compact) and $2^{\aleph_0} \leq [G^*:G^*_L] = |X^*/E_L| \leq 2^{|T|}$.
\item There is $n<\omega$ such that for every $x\in X^*,\ \diam(x) = n$. Then $E_L|_{X^*} = E_{KP}|_{X^*} = \Theta^n|_{X^*}$ and $G^*_L = X_{\Theta}^n = G'_{KP}$ are $\bigwedge$-definable groups. 
\end{enumerate}
\end{remark}

\begin{lemma} \label{lem:ind}
\begin{enumerate}
\item Either $G'_{KP} = G^*_L$, or $[G'_{KP}:G^*_L] \geq 2^{\aleph_0}$.

\item If the language of the structure $M$ is countable, then either \[ G'_{Sh} = G'_{KP} \textrm{ or } [G'_{Sh}:G'_{KP}] \geq 2^{\aleph_0}.\] In the last case the space of $\emptyset$-types $S_G(\emptyset)$ of $G$ is of power $2^{\aleph_0}$.
\end{enumerate}
\end{lemma}
\begin{proof} (1) follows from preceding remark, Lemma \ref{lem:el} and \cite[Theorem 1.1]{newelski}. 

(2) The proof is very similar to the proof of \cite[Theorem 3.5]{krupnew}, so we are brief. Consider the group $H = G^*/G'_{KP}$. This group with the logic topology is a compact Hausdorff topological group. Since the language is countable, $H$ is metrizable. Let $d_0$ be a metric on $H$. Modifying $d_0$ as in \cite{krupnew} we obtain an equivalent metric $d$, which is $\emptyset$-invariant. Since $H$ is Hausdorff, the connected component of $H$ is equal to the quasi-connected component $\QC$, and by Proposition \ref{prop:2}(5) \[\QC = G'_{Sh}/G'_{KP}.\] Assume that $G'_{Sh} \neq G'_{KP}$ and take $g\in G'_{Sh}\setminus G'_{KP}$. Let $r=d(e/G'_{KP},g/G'_{KP})$. For every $\delta$ with $0<\delta<r$ there is $g_{\delta}\in G'_{Sh}$ such that $d(e/G'_{KP},g_{\delta}/G'_{KP}) = \delta$ (because $G'_{Sh}/G'_{KP}$ is connected). The metric $d$ is $\emptyset$-invariant, hence for $\delta < \delta'$, \[\tp(g_{\delta}) \neq \tp(g_{\delta'}) \textrm{ and } d(g_{\delta'}/G'_{KP}, g_{\delta}/G'_{KP}) \geq \delta' - \delta > 0.\] Therefore the power of $S_G(\emptyset)$ is $2^{\aleph_0}$ and $g_{\delta'}g_{\delta}^{-1} \notin G'_{KP}$, hence $[G'_{Sh} : G'_{KP}] = 2^{\aleph_0}$.
\end{proof}

\section{$\bigwedge$-definability in $G$}

In this section we investigate $\bigwedge$-definability of $G^*_L$ in several special cases.

\begin{proposition}
If the theory of $M$ is small, then $G^*_L = G'_{KP} = G'_{Sh}$. Hence $G^*_L$ is $\bigwedge$-definable.
\end{proposition}
\begin{proof} Equality $G^*_L = G'_{KP}$ follows from \cite[Theorem 3.1(2)]{newelski}. Equality $G'_{KP} = G'_{Sh}$ follows from Lemma \ref{lem:ind}.
\end{proof}

\begin{proposition} \label{prop:simple}
If the theory of $M$ is simple, then the theory of $N$ is also simple and $G^*_L = X_{\Theta}^2 = G'_{KP}$.
\end{proposition}
\begin{proof} If $\Th(M)$ is simple, then $\Th(N)$ is also simple, because the structure $N' = (M,G,\cdot)$ (where $\cdot\colon G\times G \rightarrow G$ is the group action) is definable in $M$. Thus $N'$ is simple, and $N$ is obtained from $N'$ by forgetting some structure. Therefore $\Th(N)$ is also simple. In every simple structure $E_L = E_{KP} = \Theta^2$, so $G^*_L = X_{\Theta}^2$ follows from Lemma \ref{lem:el}.
\end{proof}

Now we give a criterion for equality $G'_{KP} = G'_{Sh}$, when the theory of $M$ is simple. If in this case $G'_{KP} \subsetneq G'_{Sh}$, then it gives us a solution of an open problem: there exist an example of a structure with simple theory and in which Kim-Pillay and Shelah strong types are different (see Lemma \ref{lem:el}). To state this criterion we need one definition. We call a subset $P\subseteq G^*$ \emph{thick} if $P$ is symmetric ($P = P^{-1}$) and there exist a natural number $n<\omega$ such that for any sequence $g_0,\ldots,g_{n-1}\in G$ there exist $i<j<n$ such that \[g_i^{-1}\cdot g_j \in P.\] When $\varphi(x,y)$ is a thick formula (see Definition \ref{def:thick}) then $X_{\varphi}$ (see Definition \ref{def:kom}) is thick set. On the other hand if $P$ is definable thick set, then the formula $\varphi_P(x,y) = x^{-1}\cdot y \in P$ is also thick and $P=X_{\varphi_P}$. It is easy to see that for every $n<\omega$ we have  \[X_{\Theta}^n = \bigcap \{X_{\varphi}^n : \varphi \in L \textrm{ is thick} \}.\]

\begin{lemma}
If $M$ has a simple theory, then $E_{KP}|_{X^*} \subsetneq E_{Sh}|_{X^*}$ (i.e. $G'_{KP} \subsetneq G'_{Sh}$) if and only if there exists a $\emptyset$-definable thick set $P$ such that \[G'_{Sh} \not\subseteq P^2,\] i.e. $P^2$ does not contain any $\emptyset$-definable subgroup of $G$ of finite index (see Proposition \ref{prop:1:2}(3)).
\end{lemma}
\begin{proof}
If every thick $P$ satisfies $G'_{Sh} \subseteq P^2$, then clearly $G'_{Sh} \subseteq X_{\Theta}^2 = G^*_L = G'_{KP}$, so $G'_{Sh} = G'_{KP}$.

If $G'_{Sh} = G'_{KP}$, then from \ref{prop:simple} we have that $G'_{Sh}=X_{\Theta}^2 = \bigcap \{P^2 : P \textrm{ is thick}\}$. Thus every thick $P$ satisfies $G'_{Sh} \subseteq P^2$.
\end{proof}

\begin{example}
There is an example of an abelian group $(G,\cdot,\ldots)$ which has a simple $\omega$-categorical theory and satisfies $X_{\Theta} \subsetneq X_{\Theta}^2 = G^*$ (Example 6.1.10 in \cite{cherlin}, private communication by E. Hrushovski). Consider a countable infinite dimensional vector space $V$ over $\F_2 = \{0,1\}$. Let $\B = \{b_i : i<\omega\}$ be its basis and $Q\colon V \rightarrow \F_2$ be the following degenerate orthogonal form with the induced scalar product $(\cdot,\cdot)$: \[Q\left(\sum_{i} \lambda_i b_i\right) = \lambda_0^2 + \lambda_1\lambda_2 + \lambda_3\lambda_4 + \ldots,\ \  (a,b) = Q(a+b) - Q(a) - Q(b),\ a,b\in V.\] $Q$ is degenerate, because its radical $K = \{v\in V : (v,\cdot) \equiv 0 \} = \{0,b_0\}$ is nontrivial. The structure $\G = (V,+,Q)$ has simple $\omega$-categorical theory. We show that \[X_{\Theta} \subseteq V \setminus \{b_0\}.\] If $\Theta(v,w)$, then $Q(v) = Q(w)$. Assume on the contrary that $v-w = b_0$, then $v=w+b_0$, so: \[ Q(w)=Q(v) = Q(w+b_0) = Q(w) + Q(b_0) + (w,b_0) = Q(w) + 1,\] and we reach a contradiction. 

There are only 4 types over $\emptyset$: $\tp(0), \tp(b_0), p(x), q(x)$, where $p, q$ are types of elements $v,w \neq 0, b_0$ with $Q(v) = 0, Q(w)=1$ respectively. The sets $X_{\Theta}, X^2_{\Theta}$ are $\emptyset$-invariant, so they must be a union of some sets described by above types. Consider $V_0 = \lin(b_0,b_k : k\geq 5) \prec V$. It is easy to see that
\[ (b_1,b_4) \underset{b_2 b_3 V_0}{\equiv} (b_1+b_3,b_4+b_2),\ \ (b_1,b_3) \underset{b_2 b_4 V_0}{\equiv} (b_1+b_4,b_3+b_2).\] Thus by Lemma \ref{lem:zieg}(ii) $b_3 = (b_1 + b_3) - b_1 \in X^2_{\Theta}$ and  $Q(b_3) = 0$. Also $b_1 \underset{V_0}{\equiv} b_1+b_4+b_3+b_2$, so $b_4+b_3+b_2 \in X^2_{\Theta}$ and $Q(b_4+b_3+b_2) = 1$. Therefore $V\setminus \{0,b_0\} \subseteq X^2_{\Theta}$ and then by Proposition \ref{prop:simple}  $X_{\Theta} \subsetneq X^2_{\Theta} = V$.
\end{example}


The next proposition gives us $\bigwedge$-definability of $G^*_L$ for some special groups definable in the $o$-minimal theories.

\begin{proposition} \label{prop:omin}
\begin{enumerate}
\item If $G$ is definably compact, definable in an $o$-minimal expansions of a real closed field, then $G^*_L = X^2_{\Theta} = G'_{KP} = G^{00}$.
\item If $(G,<,+,\ldots)$ is an $o$-minimal expansion of an ordered group $(G,<,+)$, then $G^* = G^*_L = X^2_{\Theta} = G^{00}$.
\end{enumerate}
\end{proposition}
\begin{proof} (1) In \cite{NIP} it is proved that under the above assumptions $G$ has {\it fsg} and there exists $G^{00}$ (the smallest definable subgroup of bounded index in $G^*$). It is also proved that $G^{00}$ is equal to \[\stab(p) = \{g\in G^* : g\cdot q = q \},\] for some (global) generic type $p(x) \in S(G^*)$. Since $p$ is a type over the model $G^*$, $\stab(p)\subseteq X^2_{\Theta}$. Therefore \[G^{00} = \stab(p) \subseteq X^2_{\Theta} \subseteq G^*_L \subseteq G'_{KP} = G^{00}_{\emptyset} = G^{00}.\]

(2) By \cite[Corollary 2.6]{petryk} we can find a global type $p(x) \in S(G^*)$, satisfying $\stab(p) = G^*$. Therefore $G^* = G^*_L = X^2_{\Theta} = G^{00}$.
\end{proof}

Case 1 from Remark \ref{rem:case} may lead us to a new example of a non-$G$-compact theory. 

There is a criterion for $\bigwedge$-definability of $G^*_L$ \cite[Theorem 3.1]{newelski}: $G^*_L$ is $\bigwedge$-definable if and only if $G^*_L = X_{\Theta}^n$ for some $n<\omega$. Thus if $X_{\Theta}$ generates a group in infinitely many steps, then $G^*_L$ is not $\bigwedge$-definable and Case 1 holds. 

We have some further partial results concerning $\bigwedge$-definability of $G^*_L$. These results involve generic subsets of $G$ and measures on $G$. They will be a part of Ph.D. thesis of the first author and appear in a forthcoming paper.

\section{More examples}

We were not able to construct an example of a group $G$, where $G_L$ is not $\bigwedge$-definable. We can try at least to construct a group $G$, where $G_{E}$ is not $\bigwedge$-definable for some equivalence relation $E$ other that $E_L$ (which gives rise to $G_L$).

It is rather easy to find such examples even in the stable case, with the relation $E$ $\bigwedge$-definable and coarser than equality of types $\equiv$. 

However even in the stable case we were not able to construct an example of $G$ where $G_{\equiv}$ is not $\bigwedge$-definable, although we conjecture such an example exists. In this case $G^*_{Sh}$ equals $G^0$, and is type definable, and equals $G_L$.

Since we are interested in finding an example where $G_L$ is not $\bigwedge$-definable, naturally we are interested in non-$\bigwedge$-definable $G_E$, where $E$ is close to $E_L$.

In this section we give only an example (Example \ref{ex:l1}), where $G_{\equiv}$ is not $\bigwedge$-definable. We could not come closer to $E_L$ than $\equiv$. We give also an example (Example \ref{ex:l2}) of a group $G$ with non-$G$-compact theory.

\begin{example} \label{ex:l1}
In \cite{macdonald} there is an example (for every $n<\omega$) of a finite group $G_n$ in which the commutators $X_{\inn(G_n)}$ generate the commutant $G'_n=[G_n,G_n]=G_{\inn(G_n)}$ in precisely $n$ steps. We expand the structure $(G_n,\cdot)$ to obtain a structure $\G_n$ satisfying \[\aut(\G_n) = \inn(G_n),\] i.e. every automorphism of $\G_n$ is an inner automorphism of $G_n$. Note that in $\G_n$ the set $X_{\equiv}$ equals $X_{\inn(G_n)}$ and generates a group in $n$ steps. Consider the product $\prod_{n<\omega} G_n$ of the groups $G_n$. We expand $\prod_{n<\omega} G_n$ to a structure $\G$ as follows. For each $k$ let $E_k$ be the equivalence relation on $\prod_{n<\omega} G_n$ given by 
\[E_k(u,v) \Leftrightarrow u(k) = v(k).\]
Then $\prod_{n<\omega} G_n/E_k$ is naturally identified with $G_k$. We expand $\prod_{n<\omega} G_n$ by the relations $E_k$, $k<\omega$, and the $\G_k$-structure on $G_k$ (identified with $\prod_{n<\omega} \G_n/E_k$). We denote the quotient map \[\prod_{n<\omega} G_n \rightarrow \prod_{n<\omega} G_n/E_k\] by $\pi_k$.

Let $\G^*$ be a large saturated extension of $\G$. We will prove that in $\G^*$, the group $\G^*_{\equiv}$ is not $\bigwedge$-definable. This boils down to proving that 
\[\pi_k[X_{\equiv}^{\G^*}] = X_{\equiv}^{\G_k}.  \leqno{(*)}\]
Indeed, suppose the above holds. Then, since $X_{\equiv}^{\G_k}$ generates a group in $\geq n$ steps, also $X_{\equiv}^{\G^*}$ generates a group in $\geq n$ steps. As $n$ is arbitrary, we get that $X_{\equiv}^{\G^*}$ generates the group $\G^*_{\equiv}$ in infinitely many steps. By \cite[Theorem 3.1(1)]{newelski}, the group $\G^*_{\equiv}$ is not $\bigwedge$-definable.

Now we prove $(*)$. $\subseteq$ is clear, since every automorphism of $\G^*$ induces an automorphism of $\G_k$. To prove $\supseteq$, consider $a,b \in \G_k$ with $b=f(a)$ for some $f\in\aut(\G_k)$. We can extend $f$ to an automorphism of $\G$ and then to $\G^*$. If $a=\pi_k(a')$ for $a'\in \G^*$, then \[b = f(a) = \pi_k(f(a')),\] and therefore $\pi_k(a'^{-1}f(a')) = a^{-1}b$.
\end{example}

Now we give an example of group $G$ whose theory is not $G$-compact, but case $2$ from Remark \ref{rem:case} holds.

\begin{example} \label{ex:l2}
First we construct a group with a large finite diameter of Lascar strong types. Let $\M_0 = (M_0, R, f)$ be a dense circular ordering (with respect to a ternary relation $R$), equipped with a function $f$, which is a cyclic bijection of $M_0$ respecting $R$, of period 3. This structure was considered in \cite{pillay} to construct the first example of a non-$G$-compact theory. Our group $M_3$ will be the disjoint union of $M_0$ and $\{0\}$ equipped with a structure of vector space over $\F_2$ (i.e. $M_3$ will be an abelian group of exponent 2) so that the addition $+$ on $M_3$ is ``independent" of $f$ and $R$.

To be more specific, let $L$ be the language consisting of a ternary relation symbol $R$, function symbols $+$ (binary) and $f$ (unary) and a constant $0$. To express ``independence" of $f$ and $+$ we define inductively a set of terms $\T$ in $L$ as follows.

\begin{definition} \label{def:l1}
Let $\T$ be the smallest set of terms of $L$ such that:
\begin{enumerate}
\item $v, f(v), f^2(v)$ are in $\T$ for every variable $v$.
\item If $\tau_1,\ldots,\tau_k$ ($k\geq 2$) are distinct terms in $\T$, then the terms $f(\tau_1+\ldots+\tau_k)$ and $f^2(\tau_1+\ldots+\tau_k)$ are in $\T$.
\end{enumerate}
Let $\T(\overline{x})$ be the set of terms in $\T$ in variables $\overline{x}$.
\end{definition}

$+$ will be interpreted as an associative operation, so we may omit parentheses in $\tau_1+\ldots+\tau_k$ in condition (2) in Definition \ref{def:l1}. $f$ will be interpreted as a cyclic function of period 3, so in Definition \ref{def:l1} there is no need to consider $f^k$ for $k\geq 3$.

\begin{definition} \label{def:l2}
Let $\CC$ be the class of $L$-structures $(V,+,0,R,f)$ such that:
\begin{enumerate}
\item $(V,+,0)$ is a vector space over $\F_2$, of infinite dimension.
\item $R$ is a circular order on the set $V^* = V \setminus \{0\}$.
\item $f$ is a cyclic bijection of $V^*$ of period 3, respecting $R$. That is, every point of $V^*$ is a cyclic point of $f$ of period 3 and $R(x,y,z)$ implies $R(f(x),f(y),f(z))$. Also, $R(x,f(x),f^2(x))$ holds in $V^*$.
\item For every $a\in V^*$, the set $\T(a) = \{\tau(a) : \tau \in \T\}$ is lineary independent.
\end{enumerate}
\end{definition}

Condition (4) in Definition \ref{def:l2} expresses the fact that in the structures in $\CC$ $f$ is ``independent" of the vector space structure.

First we show that the class $\CC$ is non-empty, in fact contains structures of arbitrary large power.

\begin{lemma} \label{lem:new}
Assume $(V,+,0)$ is a vector space over $\F_2$ and $f\colon V \rightarrow V$ satisfies condition 3. and 4. from Definition \ref{def:l2}, except for the part regarding $R$. Then the structure $(V,+,0,f)$ may be expanded to a structure in $\CC$.
\end{lemma}

\begin{proof}
Let $S^1 = \{z\in \CCC : |z|=1\}$ and let $f\colon S^1 \rightarrow S^1$ be the function defined by $f(z) = \zeta z$, where $\zeta = e^{2\pi i/3}$. Let $R$ be the anti-clockwise circular order on $S^1$. Let $(S^*, R^*, f^*)$ be a sufficently saturated extension of the structure $(S^1,R,f)$. We shall define an embedding $\Phi\colon V^* \rightarrow S^*$ such that for every $x\in V^*$, $\Phi(f(x)) = f^*(\Phi(x))$.

Let $V_0 \subseteq V^*$ be a set of representatives of the cycles of $f$. For $x\in V_0$ define $\Phi(x) \in S^*$ so that for no distinct $x,y\in V_0$ the elements $\Phi(x)$ and $\Phi(y)$ are in the same cycle of $f^*$. This is possible by the saturation of $S^*$. Extend $\Phi$ to $V^*$ by putting $\Phi(f(x)) = f^*(\Phi(x))$ and $\Phi(f^2(x)) = (f^*)^2(\Phi(x))$. In this way we have defined $\Phi$ as required.

$\Phi$ induces on $V^*$ a circular order $R$ such that the structure $(V,+,0,R,f)$ is in $\CC$.
\end{proof}

\begin{corollary} The class $\CC$ contains structures of arbitrary large power.
\end{corollary}
\begin{proof} Let $\V$ be the variety of algebrais $(V,+,0,f)$ over the language $\{+,0,f\}$ such that $(V,+,0)$ is a vector space over $\F_2$, $f(0)=0$ and $f^3 = \id$. Any free algebra in $\V$ satisfies the assumptions of Lemma \ref{lem:new}, hence may be expanded to a structure in $\CC$.
\end{proof}

\begin{lemma} $\CC$ is an elementary class with the joint embedding and amalgamation properties.
\end{lemma}
\begin{proof} The elementarity of $\CC$ is evident. We will prove that $\CC$ has the amalgamation property. The case of joint embedding property is easier. So assume $V_1, V_2$ are structures in $\CC$, with a common substructure $V_0$. We want to amalgamate them over $V_0$.

We can assume that $V_1, V_2$ are both subspaces of a vector space $V_3 = V_1 \oplus_{V_0} V_2$. Let $f' = f^{V_1} \cup f^{V_2}$. So $f'$ is a partial function on $V_3$, satisfying partially condition 3. and 4. from Definition  \ref{def:l2}. Take a large vector space $V_4$ such that $V_3$ is a subspace of $V_4$. We will find a subspace $V$ of $V_4$ containing $V_3$ and a function $f\colon V \rightarrow V$ extending $f'$ such that the structure $(V,+,0,f)$ satisfies the assumption of Lemma  \ref{lem:new}.

We define increasing sequence of subspaces $W_n \subseteq V_4$ and fnctions $f_n\colon W_n \cup f_n[W_n] \cup f^2_n[W_n] \rightarrow W_{n+1}$, $n<\omega$, such that:
\begin{enumerate}
\item $W_0 = V_3$,
\item $f' \subseteq f_0$,
\item $W_{n+1}$ is the linear span of $W_n \cup f_n[W_n] \cup f_n^2[W_n]$,
\item the set $\{f^i_0(x): x\in W_0 \setminus (V_1 \cup V_2), i\in \{1,2\}\}$ is linearly independent over $V_3$,
\item the set $\{f^i_{n+1}(x): x\in W_{n+1} \setminus (W_n \cup f_n[W_n] \cup f^2_n[W_n]), i\in\{1,2\}\}$ is linearly independent over $W_{n+1}$,
\item for $x\in W_n$, $f^3_n(x) = x$.
\end{enumerate}
The construction is straightforward. Let $V = \bigcup_n W_n$ and $f = \bigcup_n f_n$. Clearly, $V$ and $f$ satisfy our demands. As in the proof of Lemma \ref{lem:new} we can expand the structure $(V,+,0,f)$ to a structure in $\CC$, so that $V_1$ and $V_2$ are both substructures of $V$.
\end{proof}

Let $T_3$ be the model completion of $\Th(\CC)$. $T_3$ has quantifier elimination and its models are the existentially closed structures in $\CC$. We describe the 1-types in $T_3$.

Let $\C$ be a monster model of $T_3$. Let $a\in \C^* = \C \setminus \{0\}$. The type of $a$ is determined by the way in wich the linear span $\lin(\T(a))$ is circularly ordered by $R$, or even by the way in which the set $\lin(\T(a)) \cap (a,f(a))$ is linearly ordered by $R$. Here for $a\neq b \in \C^*$ \[(a,b) = \{ c\in \C^* : R(a,c,b) \},\ \ [a,b) = \{a\} \cup (a,b). \] $R$ induces on $(a,b)$ a linear ordering that we denote by $<$. By the proof of Lemma \ref{lem:new} we see that are $2^{\aleph_0}$ possible linear orderings of the set $\lin(\T(a)) \cap (a,f(a))$. So there are $2^{\aleph_0}$ complete 1-types over $\emptyset$ in $T_3$.

We say a few words about indiscernible sequences in $\C$. First, if $(a_n)_{n<\omega}$ is an infinite indiscernible sequence in $\C$, then $a_1 \in (a_0,f(a_0))$ or $a_1 \in (f^2(a_0),a_0)$. 

Secondly, we point how to construct an indiscernible sequence in $\C$. Assume $p(x) = \tp(a)$ for some $a\in \C^*$. Let $C^-(a), C^+(a)$ be a Dedekind cut in the set $\lin(\T(a)) \cap (a,f(a))$. That is, $C^-(a) < C^+(a)$ and $C^-(a) \cup C^+(a) = \lin(\T(a)) \cap (a,f(a))$.

It follows that for every $a'$ realising $p$, the corresponding sets $C^-(a')$, 
$C^+(a')$ are a Dedekind cut in the set $\lin(\T(a')) \cap (a',f(a'))$ and also the sets $f(C^-(a'))$, $f(C^+(a'))$ and $f^2(C^-(a'))$, $f^2(C^+(a'))$ are Dedekind cuts in the sets $\lin(\T(a')) \cap (f(a'),f^2(a'))$ and $\lin(\T(a')) \cap (f^2(a'),a')$, respectively.

We can find a sequence $(a_n)_{n<\omega}$ of elements of $\C^*$ such that for every $n>m$, \[ [a_n,f(a_n)) \subseteq (a_m,f(a_m)) \textrm{  and  } C^-(a_m) < [a_n,f(a_n)) < C^+(a_m). \]
Using the Ramsey theorem we can find such a sequence that is moreover indiscernible. Using indiscernible sequences like that we see that $p(x)$ is a strong Lascar type of diameter at least 3 and at most 6.

Similarly, replacing period 3 by period $n$ ($n\geq 3$), we construct a $G$-compact group $\M_n = (M_n,+,0,R_n,f_n)$ with the diameter of Lascar strong types $\geq n$.

Now let $\M = \prod_{3\leq n < \omega} \M_n$ be the product of the groups $\M_n$. Let $E_n$ be the equivalence relation on $\M$ given by $E_n(x,y) \Leftrightarrow x(n) = y(n)$. Then $\M/E_n$ is naturally identified with $\M_n$.

We consider $\M$ as a group expanded by the relations $E_n$ and the relations $R_n$ and functions $f_n$ on $\M/E_n$. Hence in $\M$ there is no finite bound on the diameter of Lascar strong types. By \cite{newelski}, $\M$ is not $G$-compact.
\end{example}

\end{document}